\documentclass[11pt]{article}
\usepackage{cite}
\usepackage{mathrsfs}
\usepackage{amsfonts}
\usepackage{amsmath}
\usepackage{amsfonts,amssymb,color}
\usepackage{dsfont}
\usepackage{curves}
\usepackage{mathrsfs}
\usepackage{pifont}
\usepackage{enumitem}
\usepackage{amssymb}
\usepackage{latexsym,amsmath,amssymb,amsfonts,epsfig,graphicx,cite,psfrag}
\usepackage{eepic,color,colordvi,amscd}
\usepackage{tikz}
\usetikzlibrary{shapes.geometric, arrows.meta, positioning, fit}
\usetikzlibrary{shapes, arrows.meta, positioning, calc}

\newtheorem{theorem}{Theorem}[section]

\newtheorem{lemma}{Lemma}[section]

\newtheorem{conjecture}{Conjecture}[section]

\newcommand{\qed}{\hfill\rule{0.5em}{0.809em}}

\def\qed{\hfill \ensuremath{\square}}

\def\pf{\noindent {\it Proof. }}

\begin{document}
\title{Sharp bounds and extremal structure of minimally strongly $k$-connected digraphs}
\author{Junran Yu$^{1}$\footnote{Email: nnuyu@foxmail.com}, \; Gregory Gutin$^{2}$\footnote{Email: g.gutin@rhul.ac.uk}, \; Hong-Jian Lai$^{3,4}$\footnote{Email: hjlai2015@hotmail.com}, \; Xiaoyan Zhang$^{1}$\footnote{Corresponding author: zhangxiaoyan@njnu.edu.cn}\\
	\small $^1$ Institute of Mathematics, School of Mathematical Sciences\\
	\small Nanjing Normal University, 1 Wenyuan Road,  Nanjing, 210023, China\\
	\small $^2$ Department of Computing,
	Security and Mathematics\\
	\small Royal Holloway University of London, London, United Kingdom\\
	\small $^3$ Department of Mathematics, West Virginia University, Morgantown, 26506, USA\\
	\small $^4$ School of Mathematics and Systems Science\\
	\small Guangdong Polytechnic Normal University, Guangzhou, 510665, China
}
\date{}
\maketitle

\begin{abstract}
	Let $D$ be a finite simple minimally strongly $k$-connected
	digraph, where $k\ge 2$, and let $n=|V(D)|$.
	Mader proved in 1985 that
	$|A(D)|\le 2k(n-k)$ for all sufficiently large $n$,
	with equality if and only if
	$D\cong \overleftrightarrow{K}_{k,n-k}$ (JCTB 1985).
	We prove that this bound holds whenever $n\ge 3k-1$
	and that the order threshold $3k-1$ is best possible.
	Moreover, if $n\ge 3k$, then equality holds if and only if
	$D\cong \overleftrightarrow{K}_{k,n-k}$.
	The threshold $3k$ is also best possible.
	These results resolve the conjecture of Mader in (JCTB 1985).
\end{abstract}

\section{Introduction}\label{sec:introduction}
All digraphs in this paper are finite and simple. A digraph $D=(V,A)$ is
\emph{loopless} if $vv\notin A$ for every $v\in V$, and it is
\emph{simple} if it is loopless and has at most one arc from $u$ to $v$ for each ordered pair of distinct vertices $u,v$. Opposite arcs $uv$ and $vu$ are allowed. We follow \cite{BM08,BG09,BG2018} for undefined notations and terminology. Let $V(D)$ and $A(D)$ for the vertex set and arc set of $D$, respectively, and call $|V(D)|$ the \emph{order} of $D$.

Minimal connectivity is a central topic in both structural and extremal graph theory. A minimally $k$-connected graph or digraph  has $k$ internally vertex-disjoint paths or directed paths between every pair of vertices in a graph and every ordered pair in a digraph, but the deletion of any edge or arc destroys this property.  The theorem obtained by Halin, which states that every finite minimally $k$-connected graph has a vertex of degree $k$, is the classical starting point for the research of existence and number of vertices of small degree in minimally $k$-connected or $k$-edge connected graphs and digraphs \cite{Halin1969}.

In the undirected setting, Mader developed Halin's theorem into a
structural and extremal theory for minimally $k$-connected graphs. In particular, he proved that the vertices of degree greater than $k$ induce a
forest \cite{Mader1972Degree,Mader1972Extremal}. Using this forest structure result, Mader proved that every minimally $k$-connected graph $G$ with $|V(G)|\ge3k-2$ satisfies $|E(G)|\le k(|V(G)|-k)$, with $K_{k,|V(G)|-k}$ as the natural extremal graph.  Later, Cai determined $|E(G)|\le \lfloor(|V(G)|+k)^2/8\rfloor$ for $k+1\leq|V(G)|\leq 3k-2$ and listed the cases when this bound is tight \cite{Mader1971,Mader1972Extremal,Cai1982}. These results show that the forest induced by the vertices of
degree greater than $k$ gives a powerful tool for obtaining
sharp extremal bounds for minimally $k$-connected graphs.

Some of these results above were extended to digraphs, especially in the work of Hamidoune \cite{Hamidoune1978,Hamidoune1981}.
For minimally strongly $k$-connected digraphs, Mader established a
directed analogue of the forest theorem \cite{Mader1985}.
Specifically, let $D_0$ be the spanning subdigraph of $D$ whose arcs are precisely those $uv\in A(D)$ for which $d_D^+(u)>k$ and
$d_D^-(v)>k$. Mader proved that the bipartite representation of $D_0$ is a forest \cite{Mader1985}.
This immediately implies that $D$ has at least $k$ vertices of outdegree $k$ and gives a linear lower bound, in terms of $|V(D)|$, on the number of vertices $v$ satisfying $d_D^+(v)=k$ or $d_D^-(v)=k$.
Mader's later work develope this obstruction, proving that every  minimally strongly $k$-connected
digraph has at least $k+1$ vertices of outdegree $k$ and at least $k+1$ vertices of indegree $k$
\cite{Mader2002}. For $k=2$, he also proved that every minimally strongly $2$-connected digraph has a vertex of indegree and outdegree $2$ \cite{Mader2002}. For $k\ge 3$, however, whether there are similar conclusions remains open \cite{Mader1996,Mader2002}. 

The directed analogue forest theorem also yields bounds on the number of arcs.  For almost every integer $n$, Mader determined the maximum number of arcs in minimally strongly $k$-connected digraph of order $n$ and characterized the extremal digraphs \cite{Mader1985}.

\begin{theorem}\cite{Mader1985}\label{thm:mader}
	Let $k\ge2$, and let $D$ be a  minimally strongly $k$-connected digraph of order $n$. Then $|A(D)|\le 2kn-k(k+1)$, with equality if and only if $D\cong\overleftrightarrow K_{k+1}$. Moreover, if $n\ge 4k+\max\{5-k,0\}$, then $|A(D)|\le 2k(n-k)$ with equality if and only if $D\cong\overleftrightarrow K_{k,n-k}$.
\end{theorem}
Mader also constructed a non-bipartite digraph of order $3k-1$
attaining $2k(n-k)$ arcs and conjectured that 

\begin{conjecture}\cite{Mader1985}\label{conj:mader}
	Every minimally strongly $k$-connected digraph $D$ of order $n\ge 3k$ satisfies $|A(D)|\le 2k(n-k)$, with equality only if
	$D\cong\overleftrightarrow{K}_{k,n-k}$.
\end{conjecture}
He also remarked that replacing his threshold in Theorem~\ref{thm:mader} by $3k$ appeared very difficult. 
A related line of work investigates how the density of a digraph
affects the number of arcs needed to preserve strong
$k$-connectivity.
Kang~\cite{Kang2019} proved that every strongly
$k$-connected digraph $D$ of order $n$ contains a strongly
$k$-connected spanning subdigraph with at most
$kn+800k(k+\overline{\Delta}(D))$ arcs, where
$\overline{\Delta}(D)$ denotes the maximum degree of the
complement of the underlying undirected graph of $D$.
In particular, every strongly $k$-connected semicomplete
digraph of order $n$ contains such a spanning subdigraph
with at most $kn+800k^2$ arcs.
These results show how the density of the digraph
can yield smaller strongly $k$-connected
spanning subdigraphs than those guaranteed by
general extremal bound in Theorem~\ref{thm:mader}. 

In this paper, we resolve the Conjecture~\ref{conj:mader}. For every $k\ge 2$, a minimally strongly $k$-connected digraph $D$ of order $n\ge 3k$ satisfies $|A(D)|\le 2k(n-k)$, with equality if and only if $D\cong\overleftrightarrow{K}_{k,n-k}$.
We further prove that the arc bound holds whenever $n\ge 3k-1$,
and that this order threshold is best possible.
In contrast, the uniqueness conclusion fails at $n=3k-1$,
hence the sharp order threshold for uniqueness is $3k$.
Our two main results are as follows.

\begin{theorem}\label{thm:sharp-size intro}
	Let $k\ge2$, and let $D$ be a minimally strongly
	$k$-connected digraph of order $n$. If $n\ge3k-1$, then
	$|A(D)|\le 2k(n-k)$.
	The order threshold $3k-1$ is best possible.
\end{theorem}

The uniqueness of the extremal digraph has the sharp order threshold $3k$.

\begin{theorem}\label{thm:sharp-equality intro}
	Let $k\ge2$, and let $D$ be a minimally strongly $k$-connected digraph of order $n\ge3k$. If
	$|A(D)|=2k(n-k)$, then $D\cong\overleftrightarrow K_{k,n-k}$.
	The order threshold $3k$ is best possible.
\end{theorem}

Both thresholds are attained by the same digraph class. For $\ell\ge2k-1$, let
$J_{k,\ell}=\overrightarrow C_{\ell}\vee\overline K_{k-1}$,
where $\vee$ adds both possible arcs between the
two vertex sets. The digraph $J_{k,\ell}$ is minimally strongly
$k$-connected and satisfies $|V(J_{k,\ell})|=\ell+k-1$ and $|A(J_{k,\ell})|=(2k-1)\ell$. Hence
$|A(J_{k,\ell})|-2k(|V(J_{k,\ell})|-k)=2k-\ell$.
Thus $J_{k,2k-1}$ has order $3k-2$ and has exactly one more arc than $2k(|V(J_{k,2k-1})|-k)$, showing that the hypothesis
$n\ge 3k-1$ in Theorem~\ref{thm:sharp-size intro} cannot be relaxed. The digraph $J_{k,2k}$, which is the equality construction recorded by Mader \cite{Mader1985},
has order $3k-1$, attains the arc bound $2k(n-k)$, and is not isomorphic to $\overleftrightarrow K_{k,2k-1}$. This proves the sharpness assertion in Theorem~\ref{thm:sharp-equality intro}.

The paper is organized as follows. Section~2 gives some notations and preliminary results used in this paper. Sections~3 and~4 prove Theorem~\ref{thm:sharp-size intro} and~\ref{thm:sharp-equality intro}, respectively. 

\section{Preliminaries}\label{sec:preliminaries}

For two vertex sets $U,W\subseteq V(D)$, define
$A_D(U,W)=\{uv\in A(D):u\in U,\ v\in W\}$.
The \emph{out-neighborhood} and \emph{in-neighborhood} of $U$ are
$N_D^+(U)=\{v\in V(D)\setminus U:
uv\in A(D)\text{ for some }u\in U\}$ and
$N_D^-(U)=\{v\in V(D)\setminus U:vu\in A(D)\text{ for some }u\in U\}$, repectively.
For a vertex $v$, we abbreviate $N_D^\pm(\{v\})$ to
$N_D^\pm(v)$ and let
$d_D^+(v)=|N_D^+(v)|$, $d_D^-(v)=|N_D^-(v)|$.
For $S\subseteq V(D)$, the digraph $D-S$ is obtained by deleting the vertices of $S$ and all arcs incident with them. For $e\in A(D)$, the digraph $D-e$ is obtained by deleting only $e$. 
For an arc $xy\in A(D)$, the vertices $x$ and $y$ are called its
\emph{tail} and \emph {head}, respectively.

A digraph is \emph{strongly connected} if it contains a directed path from $u$ to $v$ for every ordered pair of distinct vertices $u,v$. For $k\ge1$, a digraph $D$ is \emph{strongly $k$-connected} if
$|V(D)|\ge k+1$ and $D-S$ is strongly connected for every
$S\subseteq V(D)$ with $|S|\le k-1$. Equivalently, by the directed version of Menger's theorem, every ordered pair of distinct vertices is joined by $k$ internally vertex-disjoint directed paths. We denote by $\kappa(D)$ the largest integer $k$ for which $D$ is strongly $k$-connected.
Let $D$ be a digraph with $\kappa(D)=k$.
An arc $e\in A(D)$ is \emph{$k$-critical} if $D-e$ is not strongly $k$-connected.
A strongly $k$-connected digraph $D$ is
\emph{minimally strongly $k$-connected} if $D-e$ is not strongly
$k$-connected for every $e\in A(D)$. Hence every arc of a minimally strongly $k$-connected digraph is $k$-critical. In particular, $d_D^+(v)\ge k$ and $d_D^-(v)\ge k$ for every $v\in V(D)$.

The \emph{reverse digraph} $\overleftarrow D$ is obtained by
replacing every arc $uv$ by $vu$. Reversing all arcs exchanges outdegree and indegree, positive and negative fragments, and all statements having a specified direction. We shall refer to this operation as \emph{reversibility}. For an undirected graph $G$, let
$\overleftrightarrow G$ denote the digraph obtained by replacing
each edge with two opposite arcs. In particular,
$\overleftrightarrow K_{a,b}$ denotes the complete
bipartite digraph.

 We use the terminology of positive and negative fragments following Mader \cite{Mader2002}.  Let $D$ be a non-complete digraph with connectivity $\kappa(D)$.  
A nonempty set $F\subseteq V(D)$ is called a \emph{positive fragment} if $|N_D^+(F)|=\kappa(D)$
and $V(D)\setminus\bigl(F\cup N_D^+(F)\bigr)\neq\varnothing$.
A \emph{negative fragment} is defined similarly, by requiring
$|N_D^-(F)|=\kappa(D)$ and $V(D)\setminus\bigl(F\cup N_D^-(F)\bigr)\neq\varnothing$.
If $F$ is a positive fragment, then
$\overline F=V(D)\setminus\bigl(F\cup N_D^+(F)\bigr)$
is its complementary negative fragment. In particular,
$N_D^-(\overline F)=N_D^+(F)$.
Let $D$ be a minimally strongly $k$-connected digraph, where $k\ge2$, and let $n=|V(D)|$.
Since $D$ is strongly $k$-connected, we have
$d_D^+(v)$, $d_D^-(v)\ge k$ for every $v\in V(D)$.  Define
$X=\{v\in V(D):d_D^+(v)=k\}$ with $R=V(D)\setminus X$,
and $Y=\{v\in V(D):d_D^-(v)=k\}$ with $C=V(D)\setminus Y$.
Hence $R=\{v\in V(D):d_D^+(v)>k\}$ and $C=\{v\in V(D):d_D^-(v)>k\}$.
We use the notation $r=|R|$, $s=|C|$ and $p=|R\cap C|$. Following Mader, we define $E_k^>(D)=\{uv\in A(D):d_D^+(u)>k,\ d_D^-(v)>k\}$, and let $D_0=(V(D),E_k^>(D))$.  
With the notation above, $E_k^>(D)=A_D(R,C)$. We call the arcs in $E_k^>(D)$ \emph{heavy arcs}. Morever, define $h=|A_D(R,C)|=|E_k^>(D)|$ and $q=|A_D(X,Y)|$. 
Let $e=xy\in E_k^{>}(D)$ and let $H=D-e$.
Since $D$ is minimally strongly $k$-connected, $H$ is not strongly
$k$-connected.
On the other hand, deleting a single arc destroys at most one member of any family of internally vertex-disjoint directed paths with the same initial and terminal vertices.
Thus, by Menger's theorem, $H$ is strongly $(k-1)$-connected, and hence $\kappa(H)=k-1$.
Choose a set $S\subseteq V(D)$ with $|S|=k-1$ such that $H-S$ is not strongly connected.
Since $D-S$ is strongly connected, neither $x$ nor $y$ belongs to $S$. Otherwise, $H-S=D-S$, a contradiction.
Moreover, there is no directed path from $x$ to $y$ in $H-S$.
Indeed, if such a path existed, then every directed path in $D-S$
using the arc $xy$ could be replaced by a directed path in $H-S$
with the same initial and terminal vertices.
Since $D-S$ is strongly connected, this implys that $H-S$ is
strongly connected, again a contradiction.
Let $V_x$ be the set of vertices reachable from $x$ in $H-S$,
including $x$ itself, and let
$V_y=V(D)\setminus(V_x\cup S)$.
Then $x\in V_x$ and $y\in V_y$, hence both sets are nonempty.
By the definition of $V_x$, there is no arc from $V_x$ to $V_y$ in $H$. Since $D$ is obtained from $H$ by adding the arc $xy$, we obtain
$V(D)=V_x\mathbin{\cup}S\mathbin{\cup}V_y$ and
$A_D(V_x,V_y)=\{xy\}$.
We also have $N_H^{+}(V_x)\subseteq S$ and $N_H^{-}(V_y)\subseteq S$.
Deleting either of these two neighbourhoods leaves both $x$ and $y$
and destroys all directed paths from $x$ to $y$ in $H$.
Consequently, each neighbourhood has cardinality at least
$\kappa(H)=k-1$.
Since $|S|=k-1$, it follows that
$N_H^{+}(V_x)=S=N_H^{-}(V_y)$.
For such a tripartition $(V_x,S,V_y)$, define $F_x^{+}(y;S)=V_x$ and $F_x^{-}(y;S)=V_y$.
Thus $F_x^{+}(y;S)$ is a positive fragment of $D-e$ containing
the tail $x$, and $F_x^{-}(y;S)$ is its complementary negative
fragment containing the head $y$.
In particular,
$N_{D-e}^{+}(F_x^{+}(y;S))=S=N_{D-e}^{-}\bigl(F_x^{-}(y;S)\bigr)$
and
$
A_D\bigl(F_x^{+}(y;S),F_x^{-}(y;S)\bigr)=\{xy\}$.
For each $xy\in E_k^{>}(D)$, fix once and for all such a related pair
with separator $S_{xy}$, and abbreviate
$F_x^{\pm}(y)=F_x^{\pm}(y;S_{xy})$.
We shall use the following fragment-disjointness property.


\begin{lemma}\cite{Mader2002}
	\label{lem:mader-disjointness}
	Let $uv,u'v\in E_k^>(D)$ be distinct arcs with the same head $v$. Then $F_u^+(v)\cap F_{u'}^+(v)=\varnothing$.
	Similarly, if $uv,uv'\in E_k^>(D)$ are distinct arcs with the same tail $u$, then $F_u^-(v)\cap F_u^-(v')=\varnothing$.
\end{lemma}

\pf
Consider first two distinct arcs
$uv,u'v\in E_k^>(D)$. Since $v$ is the common head of two arcs in $E_k^>(D)$, we have
$d_D^-(v)>k$.
By definition, $F_u^+(v)$ and $F_{u'}^+(v)$ are the positive
fragments containing the respective tails $u$ and $u'$.
Hence \cite[Lemma~0(b)]{Mader2002} gives
$F_u^+(v)\cap F_{u'}^+(v)=\varnothing$.

For the second statement, let $uv,uv'\in E_k^>(D)$
be distinct arcs with common tail $u$, and let $\overleftarrow D$ denote the digraph obtained from $D$ by reversing every arc. The arcs $vu$ and $v'u$ then belong to $E_k^>(\overleftarrow D)$ and have the common head $u$.
For the fixed related fragment decomposition associated with $uv$,
$V(D)=F_u^+(v)\mathbin{\cup}S_{uv}\mathbin{\cup}F_u^-(v)$,
the same three sets, viewed in $\overleftarrow D-vu$, give
$V(\overleftarrow D)=F_u^-(v)\mathbin{\cup}S_{uv}\mathbin{\cup}F_u^+(v)$.
Thus $F_u^-(v)$ is the positive, or tail-side, fragment associated with the reversed arc $vu$.  Similarly, $F_u^-(v')$ is the positive fragment associated with $v'u$.
Applying the first part in $\overleftarrow D$ to the two arcs
$vu$ and $v'u$, which have the same head $u$, yields
$F_u^-(v)\cap F_u^-(v')=\varnothing$.
\qed

We next define the bipartite representation of a digraph. For each $v\in V(D)$, let $v^+=(v,+)$ and $v^-=(v,-)$
be two copies of $v$, and let $V^+=\{v^+:v\in V\}$ and $V^-=\{v^-:v\in V\}$. The bipartite representation of a digraph $D$ is defined as $B(D)=(V^+,V^-,E)$ such that $u^+v^-\in E(B(D))$ if and only if $uv\in A(D)$. Recall that $R=\{v:d_D^+(v)>k\}$ and $C=\{v:d_D^-(v)>k\}$, let
$R^+=R\times\{+\}$ and $C^-=C\times\{-\}$. Define bipartite graph $\mathcal F_D$ with parts $R^+$ and $C^-$ as follows. 
The bipartite graph $\mathcal F_D$ has $V(\mathcal F_D)=R^+\mathbin{\cup}C^-$ with $u^+v^-\in E(\mathcal F_D)$ if and only if $uv\in E_k^>(D)$. In other words, each arc in $E_k^>(D)$ is represented by an edge joining
the copy of its tail in $R^+$ to the copy of its head in $C^-$.
In particular, $|V(\mathcal F_D)|=r+s$ and $|E(\mathcal F_D)|=|E_k^>(D)|=h$. Since $D_0=(V(D),E_k^>(D))$, by definition,
$\mathcal F_D=B(D_0)[R^+\cup C^-]$.
Every vertex of $B(D_0)$ outside $R^+\cup C^-$ is isolated.
The graph $\mathcal F_D$ is precisely the heavy-arc bipartite
graph considered by Mader in \cite{Mader1985}. His forest theorem yields the following lemma.


\begin{lemma}\cite{Mader1985,Mader2002}
	\label{lem:heavy-forest}
	The bipartite graph $\mathcal F_D=B(D_0)[R^+\cup C^-]$ is a forest.  Consequently,
	if $h>0$, then
	\begin{equation}
		h\le r+s-1.
		\label{eq:heavy-forest-bound}
	\end{equation}
\end{lemma}

\pf
Mader proved that the bipartite representation $B(D_0)$ associated with
the arcs of $E_k^>(D)$ is a forest
\cite[Theorem~0(a)]{Mader2002}.
Since $\mathcal F_D=B(D_0)[R^+\cup C^-]$, the first assertion follows.
If $h>0$, then $\mathcal F_D$ is a nonempty forest on
$r+s$ vertices, and hence
$h=|E(\mathcal F_D)|\le |V(\mathcal F_D)|-1=r+s-1$.
\qed

Recall that $X\mathbin{\cup}R=V(D)$, $C\mathbin{\cup}Y=V(D)$
and that $E_k^>(D)=A_D(R,C)$ with $h=|E_k^>(D)|$, $q=|A_D(X,Y)|$.
Hence $A_D(X,C)\cup A_D(X,Y)\cup E_k^>(D)\cup A_D(R,Y)=A(D)$. The four pairs of tail and head classes $X,R$ and $C,Y$ give an
exact decomposition of the arc set.

\begin{lemma}
	\label{lem:master}
	We have $|A(D)|=2kn-k(r+s)+h-q$.
\end{lemma}

\pf
	Since every vertex of $X$ has outdegree $k$, we obtain
	$|A_D(X,C)|+q=k|X|=k(n-r)$. Similarly,
	since every vertex of $Y$ has indegree $k$, we have
	$|A_D(R,Y)|+q=k|Y|=k(n-s)$.
	Now
	$A(D)=A_D(X,C)\mathbin{\cup}A_D(X,Y)
	\mathbin{\cup}A_D(R,C)
	\mathbin{\cup}A_D(R,Y)$, and $|A(D)|=2kn-k(r+s)+h-q$ follows.
\qed

Define the two sets of missing admissible ordered pairs by
$$M^+=\bigl\{(x,c)\in X\times C:x\neq c,\ xc\notin A(D)\bigr\},$$ and
$$M^-=\bigl\{(r',y)\in R\times Y: r'\neq y,\ r'y\notin A(D)\bigr\}.$$
Let $\delta_+=|M^+|$ and $\delta_-=|M^-|$.
To compute $\delta_+$, note that the number of admissible
ordered pairs $(x,c)\in X\times C$ is
$|X||C|-|X\cap C|$.
Since
$|X|=n-r$, $|C|=s$ and $X\cap C=C\setminus R$.
we have $|X\cap C|=|C\setminus R|=s-p$. Thus the number of admissible ordered pairs in $X\times C$ is
$(n-r)s-(s-p)$.
On the other hand, we know that $|A_D(X,C)|=k(n-r)-q$.
Consequently,
\begin{equation}
	\delta_+=q+(n-r)(s-k)-(s-p)\ge0.
	\label{eq:defect-plus}
\end{equation}

Similarly, $|R\cap Y|=|R\setminus C|=r-p$, thus the number of
admissible ordered pairs in $R\times Y$ is
$r(n-s)-(r-p)$. Since $|A_D(R,Y)|=k(n-s)-q$, we obtain
\begin{equation}
	\delta_-=q+(n-s)(r-k)-(r-p)\ge0.
	\label{eq:defect-minus}
\end{equation}

\subsection{The special fragment lemma for $k=2$}
\label{subsec:k2-fragments}

For $k=2$, Mader's fragment theory gives a stronger conclusion than is available for general $k$.  Namely, every arc of $E_2^>(D)$ forces a missing admissible pair on each side.

\begin{lemma}\cite{Mader2002}
	\label{lem:mader-k2-fragment}
	Let $D$ be minimally strongly $2$-connected, and let
	$e=uv\in E_2^{>}(D)$. Then $F_u^-(v)\cap Y\neq\varnothing$
	and $F_u^+(v)\cap X\neq\varnothing$.

\end{lemma}

\pf
	The first assertion is \cite[Proposition~1]{Mader2002}, applied to the
	head-side fragment associated with $uv$. The second follows by applying the first assertion to the reversed digraph.
\qed

The result used later in the extremal proof is the following immediate consequence.

\begin{lemma}
	\label{lem:k2-special}
	If $D$ is minimally strongly $2$-connected, then
    $h\le\delta_+$ and $h\le\delta_-$.
\end{lemma}

\pf
Let $uv\in E_2^>(D)$.
Recall that $F_u^-(v)$ is the negative fragment containing the head $v$.
By Lemma~\ref{lem:mader-k2-fragment},
$F_u^-(v)\cap Y\neq\varnothing$.
Choose $y_{uv}\in F_u^-(v)\cap Y$.
Since $u\in F_u^+(v)$ and the two fragments $F_u^+(v)$, $F_u^-(v)$ are disjoint, we have $y_{uv}\neq u$.  Moreover, $v\in C$ whereas
$y_{uv}\in Y$, hence $y_{uv}\neq v$.
By the unique-crossing property, $A_D\bigl(F_u^+(v),F_u^-(v)\bigr)=\{uv\}$
and hence $uy_{uv}\notin A(D)$.
Therefore $(u,y_{uv})\in M^-$.

For a fixed tail $u$ with arc $uv, uv'\in E_2^>(D)$, the negative fragments $F_u^-(v)$, $F_u^-(v')$
corresponding to distinct arcs $uv, uv'\in E_2^>(D)$ are pairwise
disjoint by Lemma~\ref{lem:mader-disjointness}.  Hence the chosen vertices $y_{uv}\in F_u^-(v)\cap Y$ and $y_{uv'}\in F_u^-(v')\cap Y$ are distinct for arcs with the same tail. It follows that arcs with different tails give different first coordinates.
Thus $uv\mapsto (u,y_{uv})$
is an injection from $E_2^>(D)$ into $M^-$.  Consequently, $h=|E_2^>(D)|\le |M^-|=\delta_-$.

Applying the same argument to the reverse digraph gives an
injection of $E_2^>(D)$ into $M^+$, and therefore
$h\le\delta_+$.
\qed

\section{Proof of Theorem~\ref{thm:sharp-size intro}}
\subsection{Preliminary lemmas}

For each heavy arc $uv\in E_k^>(D)$, define
$\alpha_{uv}=|F_u^+(v)\cap X|$ and $\eta_{uv}=|F_u^-(v)\cap Y|$ and let $w_{uv}=\alpha_{uv}+\eta_{uv}$. We call $w_{uv}$ the \emph{charge} of the arc $uv$.
A heavy arc $uv$ is called a \emph{zero-charge arc} if
$w_{uv}=0$, and \emph{charge-one arc} if $w_{uv}=1$.
The following lemma bounds the total charge by
$\delta^++\delta^-$ and gives a lower bound on the charge of each
heavy arc.

\begin{lemma}\label{lem:general-charging}
	For every arc $uv\in E_k^>(D)$, we have 	
	$$\sum_{uv\in E_k^>(D)}\alpha_{uv}\le\delta_+  \quad \text{and}  \sum_{uv\in E_k^>(D)}\eta_{uv}\le\delta_-.$$
	Consequently, $$\sum_{uv\in E_k^>(D)}w_{uv}\le\delta_++\delta_-.$$
	Moreover, for every $uv\in E_k^>(D)$, $w_{uv}\ge n-k+1-r-s+p$. Thus, if $\lambda=n-k+1-r-s+p>0$, then
	\begin{equation}
	\lambda h\le\delta_++\delta_-.
	\label{eq:global-charging}
	\end{equation}
\end{lemma}

\pf
	Let $uv\in E_k^>(D)$ and
	$x\in F_u^+(v)\cap X$.
	Since $v\in F_u^-(v)$ and
	$A_D(F_u^+(v),F_u^-(v))=\{uv\}$,
	the arc $xv\notin A(D)$.  Moreover,
	$x\neq v$, since the two fragments
	$F_u^+(v)$ and $F_u^-(v)$ are disjoint and $x\neq u$ since $x\in X$ whereas $u\in R$.
	Thus $(x,v)\in M^+$.
	
	We associate with each choice of a heavy arc
	$uv\in E_k^{>}(D)$ and a vertex
	$x\in F_u^{+}(v)\cap X$ the missing admissible
	ordered pair $(x,v)\in M^{+}$.
	We claim that this association is injective. Indeed, the map
	$(uv,x)\mapsto (x,v)$ with $uv\in E_k^{>}(D)$ and $x\in F_u^{+}(v)\cap X$ is injective into $M^{+}$.
	Suppose that two choices $(uv,x)$ and
	$(u'v',x')$ yield the same ordered pair.
	Then $x=x'$ and $v=v'$.
	If $u\neq u'$, then $x\in F_u^{+}(v)\cap F_{u'}^{+}(v)$,
	contrary to Lemma~\ref{lem:mader-disjointness}.
	Thus $u=u'$, hence the two choices are the same.
	Since there are exactly $\alpha_{uv}$ choices of $x$
	for each heavy arc $uv$, it follows that
	$$
	\sum_{uv\in E_k^{>}(D)}\alpha_{uv}
	\le |M^{+}|=\delta_{+}.$$
	
	For the second inequality, let $uv\in E_k^{>}(D)$ and
	$y\in F_u^{-}(v)\cap Y$.
	Since $u\in F_u^{+}(v)$ and the two fragments $F_u^{+}(v)$, $F_u^{-}(v)$ are disjoint, we have $y\neq u$.
	Moreover, $y\neq v$, since $y\in Y$ whereas $v\in C$.
	It follows from
	$A_D(F_u^{+}(v),F_u^{-}(v))=\{uv\}$
	that $uy\notin A(D)$.
	Together with $u\in R$, we conclude that $(u,y)\in M^{-}$.
	
	We associate with each choice of a heavy arc
	$uv\in E_k^{>}(D)$ and a vertex
	$y\in F_u^{-}(v)\cap Y$ the missing admissible
	ordered pair $(u,y)\in M^{-}$.
	We claim that this association is injective. Indeed, the map
	$(uv,y)\mapsto (u,y)$ with $uv\in E_k^{>}(D)$ and $y\in F_u^{-}(v)\cap Y$ is injective into $M^{-}$.
	Suppose that two choices $(uv,y)$ and
	$(u'v',y')$ yield the same ordered pair.
	Then $u=u'$ and $y=y'$.
	If $v\neq v'$, then
	$y\in F_u^{-}(v)\cap F_u^{-}(v')$,
	contrary to Lemma~\ref{lem:mader-disjointness},
	which asserts that the negative fragments associated
	with distinct heavy arcs having the same tail are
	pairwise disjoint.
	Thus $v=v'$, hence the two choices are the same.
	Since there are exactly $\eta_{uv}$ choices of $y$
	for each heavy arc $uv$, it follows that
	$$
	\sum_{uv\in E_k^{>}(D)}\eta_{uv}
	\le |M^{-}|=\delta_{-}.$$
	
	Combining the two inequalities and recalling that
	$w_{uv}=\alpha_{uv}+\eta_{uv}$, we obtain
	$$
	\sum_{uv\in E_k^{>}(D)}w_{uv}
	=
	\sum_{uv\in E_k^{>}(D)}\alpha_{uv}
	+
	\sum_{uv\in E_k^{>}(D)}\eta_{uv}
	\le \delta_{+}+\delta_{-}.$$
	This proves the first statement.

	For the second statement, we first define
	$\beta_{uv}=|F_u^+(v)\cap R|$ and
	$\gamma_{uv}=|F_u^-(v)\cap C|$.
	Since $X\mathbin{\cup}R=V(D)$ and $C\mathbin{\cup}Y=V(D)$,
	we have $|F_u^+(v)|=\alpha_{uv}+\beta_{uv}$ and
	$|F_u^-(v)|=\gamma_{uv}+\eta_{uv}$.
	Using $|F_u^+(v)|+|F_u^-(v)|=n-k+1$,
	we obtain
	$\alpha_{uv}+\beta_{uv}+\gamma_{uv}+\eta_{uv}=n-k+1$.
	Since $F_u^+(v)$ and $F_u^-(v)$ are disjoint, a vertex of
	$R\cap C$ cannot be counted simultaneously in
	$\beta_{uv}$ and $\gamma_{uv}$.  Hence
	$\beta_{uv}+\gamma_{uv}\le |R\cup C|=r+s-p$.
	It follows that $$w_{uv}=\alpha_{uv}+\eta_{uv}=n-k+1-(\beta_{uv}+\gamma_{uv})\ge n-k+1-r-s+p.$$
	
	Finally, if $\lambda=n-k+1-r-s+p>0$,
	then every arc $uv\in E_k^>(D)$ satisfies
	$w_{uv}\ge\lambda$.  Since
	$|E_k^>(D)|=h$, summing over all such arcs yields
	$$\lambda h\le\sum_{uv\in E_k^>(D)}w_{uv}\le\delta_++\delta_-,$$
	as required.
\qed

\subsection{The sharp size theorem}

\begin{theorem}
	\label{thm:sharp-size}
	Let $k\ge2$, and let $D$ be a minimally strongly
	$k$-connected digraph of order $n$.
	If $n\ge3k-1$, then $|A(D)|\le2k(n-k)$.
	The order threshold $3k-1$ is best possible.
\end{theorem}

\pf
Assume first that $k\ge3$.
Suppose, for a contradiction, that
$|A(D)|\ge2k(n-k)+1$.
By Lemma~\ref{lem:master},
$|A(D)|=2kn-k(r+s)+h-q$, and hence
\begin{equation}
	h-q\ge k(r+s-2k)+1.
	\label{eq:strict-master}
\end{equation}

If $h=0$, then \eqref{eq:strict-master}, together with
$q\ge0$, gives $r+s\le2k-1$.
If $h>0$, then Lemma~\ref{lem:heavy-forest} gives
$h\le r+s-1$, hence
$k(r+s-2k)+1\le h-q\le h\le r+s-1$.
Consequently,
$(k-1)(r+s-2k)\le2k-2$.
Thus, in either case,
\begin{equation}
	r+s\le2k+2.
	\label{eq:strict-small-core}
\end{equation}
Let $n=3k-2+d$, where $d\ge1$.

\medskip
\noindent
\textbf{Case 1: both $r<k$ and $s<k$.}

Without loss of generality, we may assume that $r\le s$.
Let $r=k-a$ and $s=k-b$, where
$1\le b\le a\le k$, and let $x=a-b\ge0$.
Let $\delta=\delta_-\ge0$.
By \eqref{eq:defect-minus},
$q=a(n-s)+(r-p)+\delta$.
Define $G=a(n-s)+(r-p)-k(a+b)$.
A direct calculation gives
$G=b(b+d-3)+x(b+d+k-3)+k-p$.
Since $p\le r=k-b-x$, we obtain
$G\ge b(b+d-2)+x(b+d+k-2)\ge kx\ge0$, where the second inequality uses $b+d\ge2$.
By \eqref{eq:strict-master},
$h\ge G+\delta+1$.

The definitions of $\delta_+$ and $\delta_-$ give
$\delta_+-\delta_-=(s-r)(n-k-1)=x(n-k-1)$.
Hence
$\delta_++\delta_-=2\delta+x(n-k-1)$.
Moreover,
$\lambda=n-k+1-r-s+p=a+b+d-1+p\ge d+1\ge2$.
Since $k\lambda-(n-k-1)\ge k(d+1)-(2k+d-3)=(k-1)(d-1)+2>0$,
it follows from $G\ge kx$ and $x\ge0$ that $\lambda G\ge x(n-k-1)$.
Also, $\lambda\delta\ge2\delta$.
Therefore, $$\lambda h\ge \lambda G+\lambda\delta+\lambda\ge x(n-k-1)+2\delta+\lambda>\delta_++\delta_-,$$
contrary to Lemma~\ref{lem:general-charging}.
Thus $r$ and $s$ cannot both be smaller than $k$.

\medskip
\noindent
\textbf{Case 2: exactly one of $r$ and $s$ is smaller than $k$.}

Without loss of generality, suppose that $r<k\le s$.
Let $r=k-a$ and $s=k+c$, where
$1\le a\le k$ and $c\ge0$.
Since $\delta_-\ge0$, then \eqref{eq:defect-minus} gives $q\ge a(n-s)+(r-p)$.
Using $p\le r$ and
$n-s=2k-2+d-c\ge2k-1-c$, we obtain $q\ge a(2k-1-c)$.
Hence \eqref{eq:strict-master} implies
$h\ge q+k(c-a)+1\ge a(2k-1-c)+k(c-a)+1=a(k-1)+c(k-a)+1>0$.

If $a=k$, then $r=0$, and consequently $h=0$,
a contradiction. Hence $1\le a\le k-1$.
By Lemma~\ref{lem:heavy-forest},
$a(k-1)+c(k-a)+1\le h\le r+s-1=2k+c-a-1$, which gives
$ak+c(k-a-1)\le2k-2$.

If $2\le a\le k-1$, then $k-a-1\ge0$. Consequently, the left-hand side of $ak+c(k-a-1)\le2k-2$ is at least $ak\ge2k$, a contradiction.
Hence $a=1$, and the inequality $ak+c(k-a-1)\le2k-2$ becomes
$c(k-2)\le k-2$.
Since $k\ge3$, it follows that $c\in\{0,1\}$.

\smallskip
\noindent
\textbf{Subcase 2.1: $c=0$.}

Here $r=k-1$ and $s=k$.
Let $\delta=\delta_-\ge0$.
By \eqref{eq:defect-minus},
$q=3k-3+d-p+\delta$.
Hence \eqref{eq:strict-master} gives
$h\ge2k-2+d-p+\delta$.
On the other hand, Lemma~\ref{lem:heavy-forest} gives
$h\le2k-2$.
Combining these inequalities, we obtain
$p\ge d+\delta$.
Furthermore, we have $\lambda=d+p\ge2$ and
$\delta_++\delta_-=2k-3+d+2\delta$.
Since $p\le k-1$, we also have
$h\ge k-1+d+\delta$.
Therefore, $$\lambda h\ge2(k-1+d+\delta)=\delta_++\delta_-+d+1>\delta_++\delta_-,$$
contrary to \eqref{eq:global-charging}.

\smallskip
\noindent
\textbf{Subcase 2.2: $c=1$.}

Here $r=k-1$ and $s=k+1$.
Again let $\delta=\delta_-\ge0$.
By \eqref{eq:defect-minus}, $q=3k-4+d-p+\delta$.
Since \eqref{eq:strict-master} gives $h\ge q+1$,
while Lemma~\ref{lem:heavy-forest} gives
$h\le r+s-1=2k-1$, we obtain $q\le2k-2$.
Consequently, $p\ge k+d+\delta-2$.
Together with $p\le r=k-1$, this yields
$d+\delta\le1$.
Since $d\ge1$ and $\delta\ge0$, we must have
$d=1$, $\delta=0$ and $p=k-1$.
It follows that
$q=2k-2$ and $h=2k-1$.
Moreover, $\lambda=k-1>0$, $\delta_-=0$ and $\delta_+=4k-4$.
Thus \eqref{eq:global-charging} would imply
$(k-1)(2k-1)\le4k-4$,
which is impossible for $k\ge3$.
This excludes the case of exactly one of $r$ and $s$ is
smaller than $k$. Together with Case~1, this shows that $r\ge k$ and $s\ge k$.

\medskip
\noindent
\textbf{Case 3: $r\ge k$ and $s\ge k$.}

Combining \eqref{eq:strict-small-core} with $r,s\ge k$,
it remains, without loss of generality, to consider
$$
(r,s)\in
\{(k,k),(k,k+1),(k,k+2),(k+1,k+1)\}.$$
Throughout this case, \eqref{eq:strict-master} gives
$h\ge q+k(r+s-2k)+1\ge q+1>0$.
Hence Lemma~\ref{lem:heavy-forest} gives
$h\le r+s-1$.

\smallskip
\noindent
\textbf{Subcase 3.1: $(r,s)=(k,k)$.}

Let $\delta=q-(k-p)\ge0$.
Then we have $q=k-p+\delta$, $\delta_+=\delta_-=\delta$ and $\lambda=d+p-1\ge0$.
By \eqref{eq:strict-master} and
Lemma~\ref{lem:heavy-forest},
$k-p+\delta+1\le h\le2k-1$.
In particular,
$\delta\le k+p-2$.

If $\lambda\ge2$, then
$\lambda h\ge2(k-p+\delta+1)>2\delta=\delta_++\delta_-$, contrary to \eqref{eq:global-charging}.
Thus $\lambda\in\{0,1\}$.
Since $\lambda=d+p-1$ and $d\ge1$, we have $p\le1$.
Consequently,
$h-2\delta\ge k-p+1-\delta\ge3-2p\ge1$.
By Lemma~\ref{lem:general-charging},
$$\sum_{xy\in E_k^>(D)}w_{xy}\le2\delta<h.$$
Since all charges are nonnegative integers, there exists
an arc $uv\in E_k^>(D)$ with $w_{uv}=0$. 
By Lemma~\ref{lem:general-charging},
$w_{uv}\ge n-k+1-r-s+p=\lambda$.
Since $w_{uv}=0$ and
$\lambda=d+p-1\ge0$, we obtain
$0=w_{uv}\ge\lambda\ge0$.
Thus $\lambda=0$, and hence $d=1$ and $p=0$.
Consequently, $n=3k-1$.
Since $w_{uv}=0$, we have
$F_u^+(v)\subseteq R$ and $F_u^-(v)\subseteq C$.
Moreover,
$|F_u^+(v)|+|F_u^-(v)|=n-k+1=2k=|R\cup C|$.
Since $R\cap C=\varnothing$, it follows that
$F_u^+(v)=R$ and $F_u^-(v)=C$.
The property of fragments gives
$E_k^>(D)=A_D(R,C)=A_D(F_u^+(v),F_u^-(v))=\{uv\}$.
Thus $h=1$, contrary to
$h\ge k+\delta+1\ge k+1$.
This excludes $(r,s)=(k,k)$.

\smallskip
\noindent
\textbf{Subcase 3.2: $(r,s)=(k,k+1)$.}

Let $\delta=\delta_-=q-(k-p)\ge0$.
Then $q=k-p+\delta$.
By \eqref{eq:strict-master} and
Lemma~\ref{lem:heavy-forest}, we obtain
$2k-p+\delta+1\le h\le2k$.
Hence
$1\le p\le k$ and
$0\le\delta\le p-1$.
Furthermore,
$\lambda=n-3k+p\ge p-1\ge0$
and
$
\delta_++\delta_-
=2\delta+n-k-1
=2\delta+2k+\lambda-p-1$.
Therefore,
$$\lambda h-(\delta_++\delta_-)\ge\lambda(2k-p+\delta+1)-(2\delta+2k+\lambda-p-1)=(\lambda-1)(2k-p)+1+(\lambda-2)\delta.$$
If $\lambda\ge2$, then the right-hand side is positive,
since $2k-p\ge k>0$ and $\delta\ge0$.
This contradicts \eqref{eq:global-charging}.
Hence $\lambda\in\{0,1\}$.

If $\lambda=0$, then $\lambda\ge p-1$ and $p\ge1$
force $p=1$.
Thus $\delta=0$, and
$\lambda=n-3k+p=0$ gives $n=3k-1$.

If $\lambda=1$, then \eqref{eq:global-charging}
and the lower bound of $\lambda h-(\delta_++\delta_-)$ imply
$0\ge\lambda h-(\delta_++\delta_-)\ge1-\delta$.
Thus $\delta\ge1$.
Since $\delta+1\le p\le\lambda+1=2$,
we obtain $p=2$ and $\delta=1$.
Again $\lambda=n-3k+p$ gives $n=3k-1$.

Consequently, $n=3k-1$, and the only remaining
possibilities are
$(p,\delta)\in\{(1,0),(2,1)\}$.
In both cases, $q=k-1$ and $h=2k$.

\smallskip
\noindent
\textbf{Subcase 3.2.1: $(p,\delta)=(1,0)$.}

Here $\delta_-=0$ and $\delta_+=2k-2$.
By Lemma~\ref{lem:general-charging},
$$
\sum_{xy\in E_k^>(D)}w_{xy}
\le2k-2
<h.$$
Hence there exists an arc $uv\in E_k^>(D)$
with $w_{uv}=0$. Let $R\cap C=\{z\}$.
The two fragments satisfy
$F_u^+(v)\subseteq R$, $F_u^-(v)\subseteq C$
and $|F_u^+(v)|+|F_u^-(v)|=2k=|R\cup C|$.
Since the fragments are disjoint, they partition
$R\cup C$.

If $z\in F_u^+(v)$, then
$F_u^+(v)=R$ and $F_u^-(v)=C\setminus\{z\}$.
By the  property of the fragments,
we have
$A_D(R,C\setminus\{z\})=A_D(F_u^+(v),F_u^-(v))=\{uv\}$.
Let $xy\in E_k^>(D)\setminus\{uv\}$.
Since $x\in R$ and $y\in C$, this equality
forces $y=z$.
Thus every heavy arc other than $uv$ has head $z$.
As $D$ is simple and loopless, there are at most
$|R|-1=k-1$ such arcs.
Therefore $h\le1+(k-1)=k$.

If $z\in F_u^-(v)$, then
$F_u^+(v)=R\setminus\{z\}$ and $F_u^-(v)=C$.
Similarly, the property of fragments gives
$A_D(R\setminus\{z\},C)
=A_D(F_u^+(v),F_u^-(v))=\{uv\}$.
For any $xy\in E_k^>(D)\setminus\{uv\}$,
we have $x\in R$ and $y\in C$, hence this equality
forces $x=z$.
Thus every heavy arc other than $uv$ has tail $z$.
There are at most $|C|-1=k$ such arcs, and hence
$h\le1+k=k+1$.
In either case, $h\le k+1<2k$, contrary to $h=2k$.

\smallskip
\noindent
\textbf{Subcase 3.2.2: $(p,\delta)=(2,1)$.}

Here $\lambda=1$ and $\delta_++\delta_-=2k$.
By Lemma~\ref{lem:general-charging}, every
$uv\in E_k^>(D)$ satisfies $w_{uv}\ge n-k+1-r-s+p=\lambda=1$.
The same lemma gives
$$\sum_{uv\in E_k^>(D)}w_{uv}
\le\delta_++\delta_-=2k.$$
Since $h=2k$, it follows that $w_{uv}=1$
for every $uv\in E_k^>(D)$.
For each such arc $uv$,
$$|F_u^+(v)\cap R|+|F_u^-(v)\cap C|=n-k+1-w_{uv}=2k-1=|R\cup C|.$$
The sets on the left are disjoint subsets of
$R\cup C$, and hence
$
(F_u^+(v)\cap R)
\cup
(F_u^-(v)\cap C)
=R\cup C$.
Let $A=R\setminus C$ and $B=C\setminus R$.
Then
$|A|=k-2\ge1$ and
$|B|=k-1\ge1$.
For every $uv\in E_k^>(D)$, the preceding equality gives
$A\subseteq F_u^+(v)$ and
$B\subseteq F_u^-(v)$.
Thus two distinct heavy arcs cannot have the same head.
Otherwise, their positive fragments would both contain
the nonempty set $A$, contrary to
Lemma~\ref{lem:mader-disjointness}.
Similarly, two distinct heavy arcs cannot have the same
tail, since their negative fragments would both contain
the nonempty set $B$.
Therefore, the $2k$ edges corresponding to the heavy arcs
form a matching in $\mathcal{F}_D$.
However, the $R^+$ side has only $k$ vertices,
hence every matching in $\mathcal{F}_D$ has size at most $|R|=k$, contradicting the existence of a matching of size $2k$.
This contradiction excludes $(r,s)=(k,k+1)$.

\smallskip
\noindent
\textbf{Subcase 3.3: $(r,s)=(k,k+2)$.}

By \eqref{eq:strict-master} and
Lemma~\ref{lem:heavy-forest},
$q+2k+1\le h\le r+s-1=2k+1$.
Since $q\ge0$, we obtain $q=0$ and $h=2k+1$.
By \eqref{eq:defect-minus}, $\delta_-=p-k$.
Since $\delta_-\ge0$ and $p\le r=k$, it follows that
$p=k$.
Consequently, $\delta_-=0$, $\lambda=n-2k-1$ and $\delta_+=2(n-k-1)$.
Since $n\ge3k-1$ and $k\ge3$, we have
$\lambda\ge k-2\ge1$.
Hence \eqref{eq:global-charging} gives
$(2k+1)\lambda\le\delta_++\delta_-=2(n-k-1)=2(k+\lambda)$.
Thus
$(2k-1)\lambda\le2k$.
Since $\lambda$ is a positive integer and $k\ge3$,
we must have $\lambda=1$.
Together with $\lambda\ge k-2$, this forces
$k=3$, and hence $n=8$.
We therefore have
$(r,s,p,q,h)=(3,5,3,0,7)$,
$\delta_-=0$,
$\delta_+=8$ and
$\lambda=1$.
By Lemma~\ref{lem:general-charging}, we have
$$\sum_{uv\in E_3^>(D)}\eta_{uv}\le\delta_-=0 \ \ \  \text{and}
\sum_{uv\in E_3^>(D)}\alpha_{uv}
\le\delta_+=8.$$
Since all $\eta_{uv}$ are nonnegative, we have
$\eta_{uv}=0$ for every $uv\in E_3^>(D)$.
The same lemma gives
$w_{uv}\ge n-k+1-r-s+p=\lambda=1$
for every $uv\in E_3^>(D)$.
Thus $\alpha_{uv}=w_{uv}\ge1$ for each such arc.

There are seven heavy arcs, and the sum of the seven
positive integers $\alpha_{uv}$ is at most $8$.
Hence at least six of these arcs satisfy
$\alpha_{uv}=1$, and therefore $w_{uv}=1$.
For every such arc $uv$, the definitions of
$\alpha_{uv}$ and $\eta_{uv}$ give
$|F_u^+(v)\cap R|+|F_u^-(v)\cap C|=n-k+1-w_{uv}=5=|R\cup C|$.
The sets on the left are disjoint subsets of
$R\cup C$, hence $(F_u^+(v)\cap R)\cup(F_u^-(v)\cap C)=R\cup C$.
In particular, $C\setminus R\subseteq F_u^-(v)$. Since $|C\setminus R|=2$, the negative fragments of these arcs all
contain the same nonempty set $C\setminus R$.
By Lemma~2.1, their tails are pairwise distinct.
However, there are at least six such arcs, whereas all their tails
belong to the three-vertex set $R$.
This contradiction excludes $(r,s)=(k,k+2)$.

\smallskip
\noindent
\textbf{Subcase 3.4: $(r,s)=(k+1,k+1)$.}

By \eqref{eq:strict-master} and
Lemma~\ref{lem:heavy-forest},
$q+2k+1\le h\le r+s-1=2k+1$.
Thus $q=0$ and $h=2k+1$.
We first show that $p\ge k$.
Suppose instead that $p<k$, and set
$S=R\cap C$ and $B=C\setminus R$.
Then $|S|=p\le k-1$ and $B\subseteq X$.
Moreover, $V(D)\setminus S=B\mathbin{\cup}Y$,
where both sets are nonempty, since
$|B|=k+1-p\ge2$ and $|Y|=n-k-1\ge2k-2>0$.
Together with $q=|A_D(X,Y)|=0$, we have
$A_{D-S}(B,Y)=A_D(B,Y)\subseteq A_D(X,Y)=\varnothing$.
Hence $D-S$ is not strongly connected,
contrary to the strong $k$-connectivity of $D$.
Therefore $p\ge k$.

Let $n=3k-1+t$, where $t\ge0$.
Then $\lambda=t+p-2\ge k-2\ge1$ and $\delta_+=\delta_-=k-1+\lambda$.
By \eqref{eq:global-charging},
$(2k+1)\lambda\le2(k-1+\lambda)$.
However, $(2k+1)\lambda-2(k-1+\lambda)=(2k-1)\lambda-2k+2\ge1$,
a contradiction.
This excludes $(r,s)=(k+1,k+1)$ and completes
the proof for $k\ge3$.

	\medskip
	\noindent
	\textbf{The case $k=2$.}
	
	Assume now that $k=2$ and $n\ge5$.
	Suppose for a contradiction that
	$|A(D)|\ge4(n-2)+1=4n-7$.
	By Lemma~\ref{lem:master},
	$|A(D)|=4n-2(r+s)+h-q$,
	and hence $h-q\ge2(r+s-4)+1$.
	Thus \eqref{eq:strict-master} remains hold for $k=2$.
	If $r\le1$, then
	$$|A(D)|=\sum_{v\in V(D)}d_D^+(v)\le2(n-r)+r(n-1)=2n+r(n-3)\le3n-3.$$
	However, $4n-7>3n-3$ for $n\ge5$, a contradiction.
	Hence $r\ge2$, and by reversibility, $s\ge2$.
	Let $r=2+a$ and $s=2+b$, where $a,b\ge0$.
	By \eqref{eq:strict-master},
	$h\ge q+2(a+b)+1>0$.
	Therefore Lemma~\ref{lem:heavy-forest} gives
	$h\le r+s-1=a+b+3$.
	Combining these inequalities, we obtain
	$q+a+b\le2$.
	
	Suppose first that $a=0$.
	Then $r=2$, hence $p\le2$, and
	\eqref{eq:defect-minus} gives
	$\delta_-=q-(2-p)\le q$.
	By \eqref{eq:strict-master} and
	Lemma~\ref{lem:k2-special},
	$q+2b+1\le h\le\delta_-\le q$,
	a contradiction.
	Thus $a\ge1$, and by reversibility, $b\ge1$.
	Since $q+a+b\le2$ and $q\ge0$, it follows that
	$a=b=1$ and $q=0$.
	Thus $r=s=3$.
	By \eqref{eq:strict-master} and
	Lemma~\ref{lem:heavy-forest}, $5\le h\le5$,
	hence $h=5$.
	
	The definitions of $\delta_+$ and $\delta_-$ give
	$\delta_+=\delta_-=n-6+p$.
	By Lemma~\ref{lem:k2-special},
	$n-6+p\ge h=5$.
	Consequently, we have $\lambda=n-7+p\ge4$ and $\delta_++\delta_-=2(\lambda+1)$.
	Applying \eqref{eq:global-charging}, we obtain
	$5\lambda\le2(\lambda+1)$, or equivalently $3\lambda\le2$, contrary to $\lambda\ge4$.
	Thus $|A(D)|\le4(n-2)$, proving the desired bound
	also for $k=2$.
	
	\medskip
	\noindent
	\textbf{Sharpness.}
	
	Fix integers $k\ge2$ and $\ell\ge2k-1$.
	Let $J_{k,\ell}=\overrightarrow C_{\ell}\vee\overline K_{k-1}$ be obtained from the disjoint union
	of a directed cycle $\overrightarrow C_\ell$ and
	an independent set $P$ of size $k-1$ by adding
	both arcs between every vertex of $P$ and every
	vertex of $\overrightarrow C_\ell$.
	
	We first verify that $J_{k,\ell}$ is strongly
	$k$-connected.
	Let $S\subseteq V(J_{k,\ell})$ with $|S|\le k-1$.
	If $P\setminus S\ne\varnothing$, then at least
	$\ell-(k-1)\ge k>0$
	cycle vertices remain.
	Since every remaining vertex of $P$ is joined in
	both directions to every remaining cycle vertex,
	$J_{k,\ell}-S$ is strongly connected.
	If $P\subseteq S$, then $|P|=k-1$ and $|S|\le k-1$
	imply that $S=P$.
	In this case, $J_{k,\ell}-S$ is the entire directed
	cycle $\overrightarrow C_\ell$, which is strongly
	connected.
	Thus $J_{k,\ell}$ is strongly $k$-connected.
	
	To verify minimality, observe that every cycle
	vertex $v_i$ satisfies
	$d^+_{J_{k,\ell}}(v_i)=d^-_{J_{k,\ell}}(v_i)=k$.
	Since $P$ is independent, every arc $xy$ of
	$J_{k,\ell}$ has at least one endvertex on the
	directed cycle.
	Consequently, $d^+_{J_{k,\ell}}(x)=k$ or
	$d^-_{J_{k,\ell}}(y)=k$.
	Deleting $xy$ therefore leaves a vertex with
	outdegree or indegree $k-1$, thus
	$J_{k,\ell}-xy$ is not strongly $k$-connected.
	Hence $J_{k,\ell}$ is minimally strongly
	$k$-connected.
	
	Finally, $|V(J_{k,\ell})|=\ell+k-1$ and $|A(J_{k,\ell})|
	=\ell+2(k-1)\ell=(2k-1)\ell$.
	Therefore, $|A(J_{k,\ell})|
		-2k(|V(J_{k,\ell})|-k)
		=(2k-1)\ell-2k(\ell-1)=2k-\ell$.
	Taking $\ell=2k-1$ and let
	$n_0=|V(J_{k,2k-1})|=3k-2$, we obtain
	$|A(J_{k,2k-1})|=2k(n_0-k)+1$.
	Thus, for every $k\ge2$, the asserted bound fails
	at order $3k-2$.
	Hence the order threshold $3k-1$ is best possible.

\section{Proof of Theorem~\ref{thm:sharp-equality intro}}

\begin{theorem}
	\label{thm:sharp-equality}
	Let $k\ge2$, and let $D$ be a minimally strongly
	$k$-connected digraph of order $n\ge3k$.
	If $|A(D)|=2k(n-k)$,
	then $D\cong\overleftrightarrow K_{k,n-k}$.
	The order threshold $3k$ is best possible.
\end{theorem}

\pf
By Lemma~\ref{lem:master},
\begin{equation}
	h-q=k(r+s-2k).
	\label{eq:equality-master}
\end{equation}
Assume first that $k\ge3$.
Let $n=3k-2+d$, where $d\ge2$.
As before, let
$\lambda=n-k+1-r-s+p$.
We first show that
\begin{equation}
	r+s\le2k+2.
	\label{eq:equality-small-core}
\end{equation}
There is nothing to prove if $r+s\le2k$.
Suppose that $r+s>2k$, and let
$m=r+s-2k\ge1$.
By \eqref{eq:equality-master} and $q\ge0$,
we have $h\ge km>0$.
Hence Lemma~\ref{lem:heavy-forest} gives
$km\le h\le r+s-1=2k+m-1$.
Consequently, $(k-1)m\le2k-1$,
and therefore $m\le2+\frac{1}{k-1}<3$.
Since $m$ is an integer, we obtain $m\le2$,
proving \eqref{eq:equality-small-core}.
We next show that $r,s\ge k$.
To this end, we exclude the following two cases.

\medskip
\noindent
\textbf{Case 1: both $r<k$ and $s<k$.}

Without loss of generality, assume that $r\le s$.
Let $r=k-a$ and $s=k-b$, where
$1\le b\le a\le k$, and let $x=a-b\ge0$.
Define $\delta=\delta_-\ge0$.
By \eqref{eq:defect-minus},
$q=a(n-s)+(r-p)+\delta$.
Define
$G=a(n-s)+(r-p)-k(a+b)$.
Since $r+s-2k=-(a+b)$,
equation \eqref{eq:equality-master} gives
$h=G+\delta$.
A direct calculation yields
$G=b(b+d-3)+x(b+d+k-3)+k-p$.
Since $p\le r=k-b-x$, we obtain
$G\ge b(b+d-2)+x(b+d+k-2)\ge kx+1$,
where the last inequality uses $b\ge1$ and $d\ge2$.
The definitions of $\delta_+$ and $\delta_-$ give
$\delta_+-\delta_-=(s-r)(n-k-1)=x(n-k-1)$.
Hence
$\delta_++\delta_-=2\delta+x(n-k-1)$.
Moreover,
$\lambda=a+b+d-1+p\ge d+1\ge3$.
Since $$k\lambda-(n-k-1)\ge k(d+1)-(2k+d-3)=(k-1)(d-1)+2>0,$$
it follows that
$\lambda h=\lambda(G+\delta)\ge\lambda(kx+1)+2\delta\ge x(n-k-1)+\lambda+2\delta>\delta_++\delta_-$,
contrary to \eqref{eq:global-charging}.
Thus $r$ and $s$ cannot both be smaller than $k$.

\medskip
\noindent
\textbf{Case 2: exactly one of $r$ and $s$
	is smaller than $k$.}

Without loss of generality, suppose that $r<k\le s$.
Let $r=k-a$ and $s=k+c$, where
$1\le a\le k$ and $c\ge0$.
By \eqref{eq:defect-minus} and $\delta_-\ge0$,
$q\ge a(n-s)+(r-p)$.
Since $p\le r$ and
$n-s=2k-2+d-c\ge2k-c$, we obtain
$q\ge a(2k-c)$.
Hence \eqref{eq:equality-master} gives
$$
	h
	=q+k(c-a)
	\ge a(2k-c)+k(c-a)
	=ak+c(k-a)
	\ge ak
	>0.
$$

If $a=k$, then $r=0$ and consequently
$h=|A_D(R,C)|=0$, a contradiction.
Thus $1\le a\le k-1$.
By Lemma~\ref{lem:heavy-forest},
$ak+c(k-a)\le h\le r+s-1=2k+c-a-1$.
Hence we obtain
$
ak+c(k-a-1)+a\le2k-1$.

If $2\le a\le k-1$, then $k-a-1\ge0$,
and the left-hand side is at least
$ak+a\ge2k+2$, a contradiction.
Hence $a=1$, and the inequality becomes $k+1+c(k-2)\le2k-1$.
Since $k\ge3$, it follows that $c\in\{0,1\}$.

\smallskip
\noindent
\textbf{Subcase 2.1: $c=0$.}

Here $r=k-1$ and $s=k$.
Define $\delta=\delta_-\ge0$.
By \eqref{eq:defect-minus},
$q=3k-3+d-p+\delta$.
Equation \eqref{eq:equality-master} therefore gives
$h=q-k=2k-3+d-p+\delta$.
Furthermore,
$\lambda=d+p\ge2$ and
$\delta_++\delta_-=2k-3+d+2\delta$ hold.
Since $h>0$ and $p\le k-1$, we have
\[
\begin{aligned}
	\lambda h-(\delta_++\delta_-)
	\ge2h-(\delta_++\delta_-)
	=2k-3+d-2p
	\ge d-1>0,
\end{aligned}
\]
contrary to \eqref{eq:global-charging}.

\smallskip
\noindent
\textbf{Subcase 2.2: $c=1$.}

Here $r=k-1$ and $s=k+1$.
Again define $\delta=\delta_-\ge0$.
By \eqref{eq:defect-minus},
$q=3k-4+d-p+\delta$.
Since $r+s=2k$, equation
\eqref{eq:equality-master} gives $h=q$.
By Lemma~\ref{lem:heavy-forest},
$
q=h\le r+s-1=2k-1$.
Consequently,
$
p\ge k+d+\delta-3$.
Together with $p\le k-1$, this yields
$d+\delta\le2$.
Since $d\ge2$ and $\delta\ge0$, we must have
$d=2$, $\delta=0$ and $p=k-1$.
It follows that
$n=3k$ and $q=h=2k-1$.
The definitions of $\lambda$, $\delta_+$ and
$\delta_-$ give $\lambda=k>0$, $\delta_-=0$ and
$\delta_+=4k-2$.
Applying \eqref{eq:global-charging}, we obtain
$k(2k-1)=\lambda h\le\delta_++\delta_-=4k-2=2(2k-1)$.
Since $2k-1>0$, this implies $k\le2$, contrary to $k\ge3$.
This excludes the case of exactly one of $r$ and $s$ is smaller than $k$. Together with Case~1, this shows that $r\geq k$ and $s\geq k$.

    \medskip
\noindent
\textbf{Case 3: $r\ge k$ and $s\ge k$.}

Without loss of generality, assume that $r\le s$.
Combining \eqref{eq:equality-small-core} with
$r,s\ge k$, it remains to consider
$(r,s)\in \{(k,k),(k,k+1),(k,k+2),(k+1,k+1)\}$.
Throughout this case, let $n=3k+t$, where
$t\ge0$ is an integer.

\smallskip
\noindent
\textbf{Subcase 3.1: $(r,s)=(k,k)$.}

By \eqref{eq:defect-minus}, we may let
$\delta=\delta_-=q-(k-p)\ge0$.
Then
$q=k-p+\delta$ and $\delta_+=\delta_-=\delta$ hold.
Equation \eqref{eq:equality-master} gives $h=q$.
Moreover, $\lambda=n-k+1-r-s+p=t+p+1\ge1$.
Hence \eqref{eq:global-charging} gives
$(t+p+1)(k-p+\delta)\le2\delta$ or equivalently $(t+p+1)(k-p)+(t+p-1)\delta\le0$.

If $t+p\ge1$, then both terms on the left are
nonnegative, since $p\le k$ and $\delta\ge0$.
Since $t+p+1>0$, it follows that $p=k$.
Then we have $t+p-1=t+k-1\ge k-1>0$,
hence the same inequality forces $\delta=0$.

If $t=p=0$, then $h=q=k-p+\delta=k+\delta>0$.
The charging inequality \eqref{eq:global-charging} gives $k+\delta\le2\delta$,
and hence $\delta\ge k$.
On the other hand, Lemma~\ref{lem:heavy-forest}
gives $k+\delta=h\le r+s-1=2k-1$,
hence $\delta\le k-1$, a contradiction.
Therefore, we obtain $p=k$ and $q=h=0$.
Since $r=s=p=k$, we have $R=C$ and $X=Y$.
Consequently, $|A_D(R,R)|=h=0$ and $|A_D(X,X)|=q=0$.
Thus there are no arcs inside either $R$ or $X$.

Every vertex $x\in X$ 
 satisfies $d_D^+(x)=d_D^-(x)=k$, and all of its in-neighbours and out-neighbours lie in the $k$-vertex set $R$, every vertex of $R$ is both an out-neighbour and an in-neighbour of $x$.
Since $D$ is simple, it follows that
$N_D^+(x)=N_D^-(x)=R$ for every $x\in X$.
Hence every vertex of $X$ is joined in both
directions to every vertex of $R$, and therefore
$D\cong\overleftrightarrow{K}_{k,n-k}$.
It remains to exclude the other three possibilities
for $(r,s)$.

\smallskip
\noindent
\textbf{Subcase 3.2: $(r,s)=(k,k+1)$.}

Let $\delta=\delta_-=q-(k-p)\ge0$.
Then $q=k-p+\delta$, and
\eqref{eq:equality-master} gives
$h=q+k=2k-p+\delta\ge k>0$.
By Lemma~\ref{lem:heavy-forest},
$2k-p+\delta=h\le r+s-1=2k$.
Thus
$
0\le\delta\le p\le k$.
Furthermore, we have
$\lambda=t+p$ and $\delta_++\delta_-=2\delta+2k+t-1$.
A direct calculation gives
\[
\begin{aligned}
	\lambda h-(\delta_++\delta_-)
	=(t+p)(2k-p+\delta)
	-(2\delta+2k+t-1)
	=(2k-p-1)(t+p-1)
	+(t+p-2)\delta.
\end{aligned}
\]
If $t+p\ge2$, then
$2k-p-1\ge k-1>0$, $t+p-1\ge1$, and
$(t+p-2)\delta\ge0$.
Hence the right-hand side is positive,
contrary to \eqref{eq:global-charging}.
Therefore $t+p\le1$.
Since $t$ and $p$ are nonnegative integers, only
$(t,p)\in\{(0,0),(1,0),(0,1)\}$ remain.

\smallskip
\textbf{Subcase 3.2.1: $t=0$ and $p=0$.}

Since $0\le\delta\le p$, we have $\delta=0$.
Hence $n=3k$, $q=k$, $h=2k$, $\delta_-=0$ and $\delta_+=2k-1$.
By Lemma~\ref{lem:general-charging},
$$
\sum_{xy\in E_k^>(D)}w_{xy}
\le\delta_++\delta_-
=2k-1
<h.$$
Since all charges are nonnegative integers, there
exists an arc $uv\in E_k^>(D)$ with $w_{uv}=0$.
By the definition of $w_{uv}$,
it follows that $F_u^+(v)\subseteq R$ and $F_u^-(v)\subseteq C$.
Moreover, we have
$|F_u^+(v)|+|F_u^-(v)|=n-k+1=2k+1=|R\cup C|$.
Since the fragments are disjoint and
$R\cap C=\varnothing$, it follows that
$F_u^+(v)=R$ and
$F_u^-(v)=C$.
By the defining property of the related fragments,
$E_k^>(D)=A_D(R,C)=A_D(F_u^+(v),F_u^-(v))=\{uv\}$.
Thus $h=1$, contrary to $h=2k$.

\smallskip
\textbf{Subcase 3.2.2: $t=1$ and $p=0$.}

Again $\delta=0$. Hence
$n=3k+1$, $q=k$, $h=2k$, $\lambda=1$ and $\delta_++\delta_-=2k$.
By Lemma~\ref{lem:general-charging}, every
$uv\in E_k^>(D)$ satisfies
$w_{uv}\ge n-k+1-r-s+p=\lambda=1$.
The same lemma gives
$$
\sum_{uv\in E_k^>(D)}w_{uv}
\le\delta_++\delta_-=2k.$$
Since $h=2k$, we obtain $w_{uv}=1$ for every
$uv\in E_k^>(D)$.
For each such arc $uv$, the definitions of the
charges give $|F_u^+(v)\cap R|+|F_u^-(v)\cap C|=n-k+1-w_{uv}=2k+1=|R\cup C|$.
The sets on the left are disjoint subsets
of $R\cup C$, and hence
$(F_u^+(v)\cap R)\cup(F_u^-(v)\cap C)=R\cup C$.
Since $R\cap C=\varnothing$, this implies
$R\subseteq F_u^+(v)$ and $C\subseteq F_u^-(v)$
for every $uv\in E_k^>(D)$.

Two distinct heavy arcs with the same head
would therefore have positive fragments both
containing the nonempty set $R$, contrary to
Lemma~\ref{lem:mader-disjointness}.
Similarly, two distinct heavy arcs cannot
have the same tail, since their negative
fragments would both contain the nonempty set $C$.
Thus the $2k$ edges corresponding to the heavy
arcs form a matching in $\mathcal F_D$.
However, the $R^+$ side has only $k$ vertices,
so every matching in $\mathcal F_D$ has size
at most $k$, a contradiction.

\smallskip
\textbf{Subcase 3.2.3: $t=0$ and $p=1$.}

Here $\delta\in\{0,1\}$ and $\lambda=1$.
Moreover,
$h=2k-1+\delta$ and $\delta_++\delta_-=2k-1+2\delta
=h+\delta$.
By Lemma~\ref{lem:general-charging}, every
$uv\in E_k^>(D)$ satisfies
$w_{uv}\ge n-k+1-r-s+p=\lambda=1$
and
$$
\sum_{uv\in E_k^>(D)}w_{uv}\le h+\delta.$$
Consequently,
$$\sum_{uv\in E_k^>(D)}(w_{uv}-1)\le\delta.$$

Let
$E_1=\{uv\in E_k^>(D):w_{uv}=1\}$.
Since every arc $xy$ outside $E_1$ has
integer charge $w_{xy}\ge2$, there are at most
$\delta$ such arcs.
Hence
$|E_1|\ge h-\delta=2k-1>k$.
For every $uv\in E_1$, we have
$w_{uv}=1=\lambda$, and therefore
$|F_u^+(v)\cap R|+|F_u^-(v)\cap C|=n-k+1-w_{uv}=2k=|R\cup C|$.
Since the sets on the left are disjoint
subsets of $R\cup C$, we obtain
$(F_u^+(v)\cap R)\cup(F_u^-(v)\cap C)=R\cup C$. In particular, for every $uv\in E_1$, we have $R\setminus C\subseteq F_u^+(v)$ and
$C\setminus R\subseteq F_u^-(v)$.
Both fixed sets are nonempty, since
$|R\setminus C|=k-1>0$ and
$|C\setminus R|=k>0$.
Thus two distinct arcs of $E_1$ cannot
have the same head. Otherwise, their positive
fragments would both contain $R\setminus C$,
contrary to Lemma~\ref{lem:mader-disjointness}.
Similarly, two distinct arcs of $E_1$
cannot have the same tail, since their negative
fragments would both contain $C\setminus R$.
Therefore, the edges corresponding to
$E_1$ form a matching in $\mathcal F_D$.
Such a matching has size at most $|R|=k$,
contradicting $|E_1|>k$.
This excludes the last possibility, and hence
$(r,s)=(k,k+1)$ is impossible.

\smallskip
\noindent
\textbf{Subcase 3.3: $(r,s)=(k,k+2)$.}

Recall that $n=3k+t$, where $t\ge0$.
Equation \eqref{eq:equality-master} gives
$h=2k+q\ge2k>0$.
Hence Lemma~\ref{lem:heavy-forest} yields
$2k+q=h\le r+s-1=2k+1$.
Since $q$ is a nonnegative integer, it follows
that $q\in\{0,1\}$.
By \eqref{eq:defect-minus},
$\delta_-=q-k+p\ge0$, hence $p\ge k-q$.
Consequently,
$$\lambda=n-k+1-r-s+p=t+p-1\ge k-q-1\ge k-2\ge1.$$
The definitions of $\delta_+$ and $\delta_-$ give
$\delta_++\delta_-=2q+2k+2t+2p-2=2(k+q+\lambda)$.
Applying \eqref{eq:global-charging}, we obtain
$(2k+q)\lambda\le2(k+q+\lambda)$ or equivalently $(2k+q-2)\lambda\le2k+2q$.
If $\lambda\ge2$, then
$(2k+q-2)\lambda\ge4k+2q-4>2k+2q$,
since $k\ge3$, a contradiction.
Therefore $\lambda=1$.

We now have $t+p=2$ and $p\ge k-q\ge k-1\ge2$.
Thus $p=2$ and $t=0$.
Since $p\ge k-q$, $q\le1$, and $k\ge3$,
we must also have $k=3$ and $q=1$.
Consequently, it follows that
$n=9$, $(r,s,p,q,h)=(3,5,2,1,7)$
with $\delta_-=0$,
$\delta_+=10$ and $\lambda=1$.
By Lemma~\ref{lem:general-charging}, we have
$$
\sum_{uv\in E_3^>(D)}\eta_{uv}\le0 \  \  \
\text{and}
\sum_{uv\in E_3^>(D)}\alpha_{uv}\le10.
$$
Since every $\eta_{uv}$ is nonnegative, hence we have
$\eta_{uv}=0$ for every $uv\in E_3^>(D)$.
The same lemma gives $w_{uv}\ge n-k+1-r-s+p=\lambda=1$.
Hence $\alpha_{uv}=w_{uv}\ge1$ for every
$uv\in E_3^>(D)$, and
$$
\sum_{uv\in E_3^>(D)}(\alpha_{uv}-1)
\le10-7=3.$$
Since all charges are integers, at most three
of the seven heavy arcs have $\alpha_{uv}\ge2$.
Therefore at least four heavy arcs satisfy
$\alpha_{uv}=w_{uv}=1$.
For every such arc $uv$, the definitions of
the charges give $|F_u^+(v)\cap R|+|F_u^-(v)\cap C| =n-k+1-w_{uv}=6=|R\cup C|$.
The sets on the left are disjoint subsets
of $R\cup C$, hence
$(F_u^+(v)\cap R)\cup(F_u^-(v)\cap C)=R\cup C$.
In particular, $C\setminus R\subseteq F_u^-(v)$
for each of these arcs.
Since $|C\setminus R|=3>0$, two such arcs
with the same tail would have negative
fragments containing the same nonempty set
$C\setminus R$, contrary to
Lemma~\ref{lem:mader-disjointness}.
Thus these arcs have pairwise distinct tails.
However, all their tails belong to $R$,
which has only three vertices.
This contradiction excludes $(r,s)=(k,k+2)$.

\smallskip
\noindent
\textbf{Subcase 3.4: $(r,s)=(k+1,k+1)$.}

By \eqref{eq:equality-master},
$h=2k+q\ge2k>0$.
Lemma~\ref{lem:heavy-forest} therefore gives
$2k+q=h\le r+s-1=2k+1$.
Hence $q\in\{0,1\}$.
Furthermore, $\lambda=t+p-1$ and $\delta_+=\delta_-=q+k+\lambda-1$ hold.
If $\lambda\ge2$, then
	$\lambda h-(\delta_++\delta_-)
	=\lambda(2k+q)-2(q+k+\lambda-1)
	=2(k-1)(\lambda-1)+q(\lambda-2)>0$,
contrary to \eqref{eq:global-charging}.
Thus $\lambda\le1$, and hence
$t+p=\lambda+1\le2$.
We next claim that $p+q\ge k$.
Suppose instead that $p+q<k$.
Let $B=C\setminus R$
and $T=\{y\in Y: by\in A(D)\text{ for some }b\in B\}$.
Since $B\subseteq X$, we have
$|T|\le |A_D(B,Y)|\le |A_D(X,Y)|=q$.
Let $S=(R\cap C)\cup T$.
Then $|S|=p+|T|\le p+q\le k-1$.
Moreover, $V(D)\setminus S=B\mathbin{\cup}(Y\setminus T)$,
where both sets are nonempty, since
$|B|=k+1-p\ge2$
and
$|Y\setminus T|\ge N-k-1-q\ge2k-2>0$.
By the definition of $T$, there is no arc
from $B$ to $Y\setminus T$ in $D-S$.
Hence $D-S$ is not strongly connected,
contrary to the strong $k$-connectivity of $D$.
This proves $p+q\ge k$.

Since $t+p\le2$, $t\ge0$ and $q\le1$, we obtain
$3\le k\le p+q\le3$. It follows that
 $k=3$, $p=2$, $q=1$ and $t=0$.
Consequently, $n=9$, $(r,s,p,q,h)=(4,4,2,1,7)$ with $\lambda=1$ and $\delta_+=\delta_-=4$.
By Lemma~\ref{lem:general-charging}, every
$uv\in E_3^>(D)$ satisfies
$w_{uv}\ge n-k+1-r-s+p=\lambda=1$
and
$$\sum_{uv\in E_3^>(D)}w_{uv}\le\delta_++\delta_-=8.$$
Since $h=7$, it follows that
$$\sum_{uv\in E_3^>(D)}(w_{uv}-1)\le1.$$
Thus at most one heavy arc has charge
greater than $1$, and at least six heavy
arcs have charge exactly $1$.
For every such arc $uv$, we have
$|F_u^+(v)\cap R|+|F_u^-(v)\cap C|=n-k+1-w_{uv}=6=|R\cup C|$.
Since the sets on the left are disjoint
subsets of $R\cup C$, their union is
$R\cup C$.
In particular,
$C\setminus R\subseteq F_u^-(v)$.
Since $|C\setminus R|=2>0$,
Lemma~\ref{lem:mader-disjointness} implies
that these arcs have pairwise
distinct tails.
This is impossible, since all their
tails lie in $R$ and $|R|=4$.
Therefore $(r,s)=(k+1,k+1)$ is impossible.
The only remaining possibility is
$(r,s)=(k,k)$.
By Subcase~3.1,
$D\cong\overleftrightarrow{K}_{k,n-k}$.
This proves the equality characterization
for $k\ge3$.

\medskip
\noindent
\textbf{The case $k=2$.}

Assume now that $k=2$ and $n\ge6$.
Recall that $|A(D)|=4(n-2)$.
Equation \eqref{eq:equality-master} becomes
$h-q=2(r+s-4)$.
If $r\le1$, then
\[
\begin{aligned}
	|A(D)|
	=\sum_{v\in V(D)}d_D^+(v)
	\le2(n-r)+r(n-1)
	=2n+r(n-3)
	\le3n-3<4n-8,
\end{aligned}
\]
where the last inequality follows from $n\ge6$.
This contradicts $|A(D)|=4(n-2)$. Hence $r\ge2$ and by reversibility, $s\ge2$.
Let $r=2+a$ and $s=2+b$, where $a,b\ge0$.
By \eqref{eq:equality-master}, $h=q+2(a+b)$.

Suppose first that $a=0$.
Then $r=2$, hence $p\le2$.
By \eqref{eq:defect-minus}, $\delta_-=q-(2-p)=q+p-2$.
Applying Lemma~\ref{lem:k2-special}, we obtain
$q+2b=h\le\delta_-=q+p-2\le q$.
Thus $2b\le p-2\le0$. Since $b\ge0$, it follows that $b=0$ and $p=2$. By reversibility, the same conclusion
$a=b=0$ and $p=2$ holds if $b=0$.

It remains to exclude the possibility that $a,b\ge1$.
Suppose that this is the case.
Then $h=q+2(a+b)\ge4>0$.
Hence Lemma~\ref{lem:heavy-forest} gives
$q+2(a+b)=h\le r+s-1=a+b+3$ and therefore $q+a+b\le3$.
Since $a,b$ are positive integers and $q\ge0$, the only
possibilities, up to reversibility, are
$(a,b)\in\{(1,1),(1,2)\}$.

If $(a,b)=(1,1)$, then $r=s=3$, $q\le1$ and $h=q+4$.
The definitions of $\delta_+$ and $\delta_-$ give
$\delta_+=\delta_-=q+n-6+p$.
By Lemma~\ref{lem:k2-special}, $q+4=h\le\delta_-=q+n-6+p$,
hence $n-6+p\ge4$.
Consequently, we obtain $\lambda=n-7+p\ge3$ and
$\delta_++\delta_-=2(q+\lambda+1)$.
However,
$$\lambda h-(\delta_++\delta_-)=\lambda(q+4)-2(q+\lambda+1)
=(\lambda-2)q+2\lambda-2\ge4>0,$$
contrary to \eqref{eq:global-charging}.

If $(a,b)=(1,2)$, then $q+a+b\le3$ gives $q=0$.
Thus $r=3$, $s=4$ and $h=6$.
The definitions of $\delta_+$ and $\delta_-$ give
$\delta_+=2n-10+p$ and $\delta_-=n-7+p$.
Moreover,
$\lambda=n-8+p$.
By Lemma~\ref{lem:k2-special}, we have $\lambda+1=\delta_-\ge h=6$,
hence $\lambda\ge5$.
We also have
$\delta_++\delta_-=3n-17+2p=3\lambda+7-p$.
Hence $\lambda h-(\delta_++\delta_-)=6\lambda-(3\lambda+7-p)=3\lambda-7+p\ge8>0$,
again contradicting \eqref{eq:global-charging}. Hence it is impossible for both $a$ and $b$ are positive.
The analysis at the begining of the cases $a=0$ and $b=0$ therefore gives $r=s=p=2$.
Thus $R=C$, $X=Y$ and $h=q$.
Furthermore, $\delta_+=\delta_-=q$ and
$\lambda=n-3\ge3$.
Applying \eqref{eq:global-charging}, we obtain
$(n-3)q\le2q$ or equivalently $(n-5)q\le0$.
Since $n\ge6$ and $q\ge0$, it follows that
$q=h=0$.
Consequently,
$|A_D(R,R)|=h=0$ and $|A_D(X,X)|=q=0$.
Thus there are no arcs inside either $R$ or $X$.
Every vertex $x\in X$ has indegree and outdegree
exactly $2$, and all of its in-neighbours and
out-neighbours lie in the two-vertex set $R$.
Since $D$ is simple, we obtain
$N_D^+(x)=N_D^-(x)=R$ for every $x\in X$.
Hence $D\cong\overleftrightarrow{K}_{2,n-2}$.
Together with the case $k\ge3$, this proves
the equality characterization for every $k\ge2$.

\medskip
\noindent
\textbf{Sharpness of the uniqueness threshold.}

Fix an integer $k\ge2$.
By the construction and verification in the proof
of Theorem~\ref{thm:sharp-size}, the digraph
$J_{k,2k}=\overrightarrow C_{2k}\vee\overline K_{k-1}$ is minimally strongly $k$-connected.
Moreover, we have $|V(J_{k,2k})|=3k-1$
and $|A(J_{k,2k})|=2k(2k-1)=2k(|V(J_{k,2k})|-k)$.
However,
$J_{k,2k}\not\cong\overleftrightarrow{K}_{k,2k-1}$.
Indeed, every vertex of the independent set
$P$ in $J_{k,2k}$ has outdegree $2k$.
Since $|P|=k-1\ge1$, such a vertex exists.
In contrast, every vertex of
$\overleftrightarrow{K}_{k,2k-1}$ has outdegree
either $k$ or $2k-1$, both of which are
smaller than $2k$.
Thus the two digraphs have different outdegree
sequences and are not isomorphic.

Therefore, for every $k\ge2$, there is a minimally
strongly $k$-connected digraph of order $3k-1$
that attains the bound but not isomorphic to complete bipartite digraph $\overleftrightarrow{K}_{k,2k-1}$.
Hence the order threshold $3k$ in the uniqueness
statement is best possible.
\qed

	{\section*{AI Declarations}
	ChatGPT was used for refining the presentation of
	this paper and translating the original papers \cite{Mader1971,Mader1972Degree,Mader1972Extremal,Mader1985,Hamidoune1981} written in German and French. All mathematical arguments and proofs in this paper were constructed and verified by the authors, who take full responsibility for all results.	
		
}

\end{document}